\documentclass[12pt,reqno]{amsart}
\usepackage{amsmath,amsthm,amssymb,amscd}

\begin{document}

\newcommand{\DATE}{March 27, 2003}
\newcommand{\TITLE}{A Lower Bound for the Canonical Height on Elliptic Curves over Abelian Extensions}
\newcommand{\TITLERUNNING}{A Lower Bound for the Canonical Height}
\newcommand{\AUTHOR}{Joseph H. Silverman}

\theoremstyle{plain}
\newtheorem{theorem}{Theorem}
\newtheorem{conjecture}[theorem]{Conjecture}
\newtheorem{proposition}[theorem]{Proposition}
\newtheorem{lemma}[theorem]{Lemma}
\newtheorem{corollary}[theorem]{Corollary}

\theoremstyle{definition}
\newtheorem*{definition}{Definition}

\theoremstyle{remark}
\newtheorem{remark}{Remark}
\newtheorem{example}{Example}
\newtheorem*{acknowledgement}{Acknowledgements}

\def\BigStrut{\vphantom{$(^{(^(}_{(}$}} 

\newenvironment{notation}[0]{%
  \begin{list}%
    {}%
    {\setlength{\itemindent}{0pt}
     \setlength{\labelwidth}{4\parindent}
     \setlength{\labelsep}{\parindent}
     \setlength{\leftmargin}{5\parindent}
     \setlength{\itemsep}{0pt}
     }%
   }%
  {\end{list}}

\newenvironment{parts}[0]{%
  \begin{list}{}%
    {\setlength{\itemindent}{0pt}
     \setlength{\labelwidth}{1.5\parindent}
     \setlength{\labelsep}{.5\parindent}
     \setlength{\leftmargin}{2\parindent}
     \setlength{\itemsep}{0pt}
     }%
   }%
  {\end{list}}
\newcommand{\Part}[1]{\item[\upshape#1]}

\renewcommand{\a}{\alpha}
\renewcommand{\b}{\beta}
\newcommand{\g}{\gamma}
\renewcommand{\d}{\delta}
\newcommand{\e}{\epsilon}
\newcommand{\f}{\phi}
\renewcommand{\l}{\lambda}
\renewcommand{\k}{\kappa}
\newcommand{\lhat}{\hat\lambda}
\newcommand{\m}{\mu}
\renewcommand{\o}{\omega}
\renewcommand{\r}{\rho}
\newcommand{\rbar}{{\bar\rho}}
\newcommand{\s}{\sigma}
\newcommand{\sbar}{{\bar\sigma}}
\renewcommand{\t}{\tau}
\newcommand{\z}{\zeta}

\newcommand{\D}{\Delta}
\newcommand{\F}{\Phi}

\newcommand{\gp}{{\mathfrak{p}}}
\newcommand{\gP}{{\mathfrak{P}}}

\def\Acal{{\mathcal A}}
\def\Bcal{{\mathcal B}}
\def\Ccal{{\mathcal C}}
\def\Dcal{{\mathcal D}}
\def\Ecal{{\mathcal E}}
\def\Fcal{{\mathcal F}}
\def\Gcal{{\mathcal G}}
\def\Hcal{{\mathcal H}}
\def\Ical{{\mathcal I}}
\def\Jcal{{\mathcal J}}
\def\Kcal{{\mathcal K}}
\def\Lcal{{\mathcal L}}
\def\Mcal{{\mathcal M}}
\def\Ncal{{\mathcal N}}
\def\Ocal{{\mathcal O}}
\def\Pcal{{\mathcal P}}
\def\Qcal{{\mathcal Q}}
\def\Rcal{{\mathcal R}}
\def\Scal{{\mathcal S}}
\def\Tcal{{\mathcal T}}
\def\Ucal{{\mathcal U}}
\def\Vcal{{\mathcal V}}
\def\Wcal{{\mathcal W}}
\def\Xcal{{\mathcal X}}
\def\Ycal{{\mathcal Y}}
\def\Zcal{{\mathcal Z}}

\renewcommand{\AA}{\mathbb{A}}
\newcommand{\BB}{\mathbb{B}}
\newcommand{\CC}{\mathbb{C}}
\newcommand{\FF}{\mathbb{F}}
\newcommand{\PP}{\mathbb{P}}
\newcommand{\QQ}{\mathbb{Q}}
\newcommand{\RR}{\mathbb{R}}
\newcommand{\ZZ}{\mathbb{Z}}

\def \bfa{{\mathbf a}}
\def \bfb{{\mathbf b}}
\def \bfc{{\mathbf c}}
\def \bfe{{\mathbf e}}
\def \bff{{\mathbf f}}
\def \bfF{{\mathbf F}}
\def \bfg{{\mathbf g}}
\def \bfp{{\mathbf p}}
\def \bfr{{\mathbf r}}
\def \bfs{{\mathbf s}}
\def \bft{{\mathbf t}}
\def \bfu{{\mathbf u}}
\def \bfv{{\mathbf v}}
\def \bfw{{\mathbf w}}
\def \bfx{{\mathbf x}}
\def \bfy{{\mathbf y}}
\def \bfz{{\mathbf z}}


\newcommand{\ab}{{\textup{ab}}}
\newcommand{\Aut}{\operatorname{Aut}}
\newcommand{\Cond}{{\mathfrak{N}}} 
\newcommand{\Disc}{\operatorname{Disc}}
\newcommand{\Div}{\operatorname{Div}}
\newcommand{\End}{\operatorname{End}}
\newcommand{\Frob}{\operatorname{Frob}}
\newcommand{\GK}{G_{\Kbar/K}}
\newcommand{\GLK}{G_{L/K}}
\newcommand{\GL}{\operatorname{GL}}
\newcommand{\Gal}{\operatorname{Gal}}
\newcommand{\GSp}{\operatorname{GSp}}
\newcommand{\hhat}{{\hat h}}
\newcommand{\Image}{\operatorname{Image}}
\newcommand{\Kbar}{{\bar K}}
\newcommand{\MOD}[1]{~(\textup{mod}~#1)}
\newcommand{\Norm}{\operatorname{N}}
\newcommand{\notdivide}{\nmid}
\newcommand{\ord}{\operatorname{ord}}
\newcommand{\Pic}{\operatorname{Pic}}
\newcommand{\Qbar}{{\bar{\QQ}}}
\newcommand{\rank}{\operatorname{rank}}
\newcommand{\res}{\operatornamewithlimits{res}}
\newcommand{\Resultant}{\operatorname{Resultant}}
\renewcommand{\setminus}{\smallsetminus}
\newcommand{\Support}{\operatorname{Support}}
\newcommand{\tors}{{\textup{tors}}}
\newcommand{\<}{\langle}
\renewcommand{\>}{\rangle}



\title[\TITLERUNNING]{\TITLE}
\date{\DATE}
\author{\AUTHOR}
\address{Mathematics Department, Box 1917, Brown University,
Providence, RI 02912 USA}
\email{jhs@math.brown.edu}
\thanks{This research supported by NSF DMS-9970382.}

\begin{abstract}
Let $E/K$ be an elliptic curve defined over a number field,
let~$\hhat$ be the canonical height on~$E$, and let~$K^\ab/K$ be the
maximal abelian extension of~$K$. Extending work of
Baker~\cite{Baker}, we prove that there is a constant~$C(E/K)>0$ so
that every nontorsion point~$P\in E(K^\ab)$
satisfies~$\hhat(P)>C(E/K)$.
\end{abstract}

\maketitle


\section{Introduction}
The classical Lehmer conjecture~\cite{Lehmer} asserts that there is an
absolute constant~$C>0$ so that any algebraic number~$\a$ that is not
a root of unity satisfies $h(\a)>C/[\QQ(\a){:}\QQ]$. Recently Amoroso
and Dvornicich~\cite{AmorosoDvornicich} showed that if~$\a$ is
restricted to lie in~$\QQ^\ab$, then the stronger inequality~$h(\a)>C$
is true.  The analog of Lehmer's conjecture for elliptic curves and
abelian varieties has also been much
studied~\cite{AndersonMasser,DavidHindry,HindrySilverman,Laurent,MasserEC,MasserAV,SilvermanDuke}. Baker~\cite{Baker}
has proven the elliptic analog of the Amo\-ro\-so-Dvor\-ni\-cich
result if the elliptic curve either has complex multiplication or has
non-integral $j$-invariant, but he was unable to handle the general
case of integral $j$-invariant. In this note we extend Baker's result
to include all elliptic curves.

\begin{theorem}
\label{theorem:maintheorem}
Let $K/\QQ$ be a number field, let~$E/K$ be an elliptic curve, 
and let~$\hhat:E(\Kbar)\to\RR$ be the canonical height on~$E$.
There is a constant $C(E/K)>0$ such that every nontorsion point
$P\in E(K^\ab)$ satisfies
$$
  \hhat(P) > C.
$$
\end{theorem}

\begin{remark}
This theorem gives a proof of Baker's Conjecture~1.10~\cite{Baker}
for elliptic curves. In Baker's terminology, we will prove
that if~$E$ does not have complex multiplication, then~$E(K^\ab)$
has the ``strong discreteness property.''
\end{remark}

We begin in Section~\ref{section:history} with a brief history of
known results on Lehmer's conjecture. In
Section~\ref{section:preliminaries} we set notation and recall basic
facts about local heights.  Section~\ref{section:torsion} discusses
torsion points and quotes a result of Serre that will be needed in the
proof of Theorem~\ref{theorem:maintheorem}. In
Section~\ref{section:unramifiedcase} we take up the unramified
case. Here we use the characteristic polynomial of
Frobenius at a prime ideal~$\gP$, evaluated at Frobenius, to annihilate points
modulo~$\gP$. This device replaces the use of a complex multiplication
map in~\cite{Baker}, allowing us to deal with general elliptic
curves. Next, in Section~\ref{section:ramifiedcase}, we deal with the
ramified case. We simplify the argument in~\cite{Baker} by
applying~\text{$\t-1$} twice to our point, the first time to move it
into the formal group, and the second time to make it a difference of
conjugate points in the formal group.  Finally, in
Section~\ref{section:proofofmaintheorem}, we complete the proof of our
main theorem. In Section~\ref{section:abelianvarieties}, we indicate
how many of the arguments can be generalized to the case of abelian
varieties. In particular, we sketch a proof that $\hhat(P)>C(A/K)$
for all nontorsion points $P\in A(K^\ab)$ if the abelian variety~$A/K$
has no complex multiplication abelian subvarieties.

\begin{acknowledgement}
The author would like to thank Matt Baker for making his
paper~\cite{Baker} available and for his many helpful comments on the
first draft of this paper, Marc Hindry for many pleasant hours spent
working together on lower bounds for canonical heights, and Michael
Rosen for numerous mathematical discussions.
\end{acknowledgement}

\section{Earlier Results on Lehmer's Conjecture}
\label{section:history}

For the convenience of the reader, in this section we summarize 
some of the known results regarding Lehmer's conjecture for number
fields and for elliptic curves.  Detailed information and an extensive
bibliography for the former is given on the Lehmer Conjecture
Website~\cite{LehmerConjHomePage}.

\subsection{The Classical Lehmer Conjecture}
Let $\a\in\Qbar^*$, assume that~$\a$ is not a root of unity, and
let~$D=[\QQ(\a){:}\QQ]$. The \emph{Lehmer Conjecture} says that there is
a constant~$c>0$ so that $h(\a)\ge cD^{-1}$ for all such~$\a$. 
\par
Blanksby and Montgomery~\cite{BlanksbyMontgomery} and
Stewart~\cite{Stewart} independently, and by rather
different methods, proved that
$$
  h(\a) \ge cD^{-1}(\log D)^{-1}.
$$
This was superceded by Dobrowolski's estimate~\cite{Dobrowolski}
$$
  h(\a) \ge cD^{-1}\bigl(\log\log(D)/\log(D)\bigr)^3,
$$
which is currently the best known general result. In the special case
that~$\a^{-1}$ is not a Galois conjugate of~$\a$ (one says that~$\a$
is non-reciprocal), Smyth~\cite{Smyth} proved the full Lehmer
conjecture~$h(\a)\ge cD^{-1}$. See also~\cite{AmorosoDvornicich} for
an alternative proof of Smyth's result.

\subsection{Lehmer's Conjecture for Elliptic Curves}
Let~$E/K$ be an elliptic curve, let~$\hhat$ be the canonical height
on~$E$, let~$P\in E(\Kbar)$ be a nontorsion point, and let
\text{$D=[K(P){:}K]$} be the degree of the field of definition of~$P$.
The elliptic analog of Lehmer's conjecture says that $\hhat(P)\ge
cD^{-1}$. Table~\ref{table:lowerbdinKbarhistory} summarizes the
history of lower bounds for~$\hhat(P)$, where~$c$ denotes a positive
constant that depends on~$E/K$, but not on~$P$.  We also note that
Pacheco~\cite[Theorem~2.3]{Pacheco} has proven the elliptic Lehmer
conjecture for function fields over finite fields.

\begin{table}
\begin{center}
\begin{tabular}{|c|c|c|} \hline
\BigStrut
  $\hhat(P)\ge{}$ & Restriction on $E$ & Reference \\ \hline\hline
\BigStrut
  $cD^{-10}(\log D)^{-6}$ & none & Anderson-Masser (1980)
       \cite{AndersonMasser} \\ \hline
\BigStrut
  $cD^{-1}\left(\frac{\log\log(D)}{\log(D)}\right)^3$ 
    & CM & Laurent (1983) \cite{Laurent} \\ \hline
\BigStrut
  $cD^{-3}(\log D)^{-2}$ & none & Masser (1989)
       \cite{MasserEC} \\ \hline
\BigStrut
  $cD^{-2}(\log D)^{-2}$ & $j$ nonintegral & 
     Hindry-Silverman (1990) \cite{HindrySilverman} \\ \hline
\end{tabular}
\caption{History of lower bounds for $\hhat$ in $E(\Kbar)$}
\label{table:lowerbdinKbarhistory}
\end{center}
\end{table}

\subsection{Lehmer's Conjecture in Abelian Extensions}
In the classical setting, Amoroso and
Dvornicich~\cite{AmorosoDvornicich} recently gave a definitive
lower bound for every $\a\in\QQ^\ab$ that is not a root of unity:
$$
  h(\a) \ge c.
$$
Amoroso and Zannier~\cite{AmorosoZannier} generalize this to relative
abelian extensions $K^\ab/K$. More generally, they give a Dobrowolski-type
bound for~$\a\in \Kbar$ with the degree $D=[K(\a){:}K]$ replaced by 
the ``nonabelian'' part of the degree $D'=[K(\a){:}K(\a)\cap K^\ab]$.
\par
Table~\ref{table:lowerbdinKabhistory} gives
the history of lower bounds for the canonical height in~$E(K^\ab)$ 
for elliptic curves~$E/K$. In particular, the
results in the present paper complete the one case left undone
in~\cite{Baker}, namely non-CM elliptic curves with integral
$j$-invariant. We observe that the non-CM integral $j$-invariant
case often presents the greatest difficulties when studying 
algebraic points on elliptic curves, see for example Serre's
theorem~\cite{SerreImageOfGalois} on the image of Galois
in~$\End(E_\tors)$. The reason is that if~$E$ has CM or if~$E$
has nonintegral $j$-invariant, then there is some control over
the action of Galois on torsion. In the former case, there is
almost complete control via class field theory, and in the
latter case, the Tate curve over~$K_\gp$ gives control for
the decomposition group at~$\gp$.

\begin{table}
\begin{center}
\begin{tabular}{|c|c|c|} \hline
\BigStrut
  $\hhat(P)\ge{}$ & Restriction on $E$ & Reference \\ \hline\hline
\BigStrut
  $cD^{-2}$ & none & Silverman (1981) \cite{SilvermanDuke} \\ \hline
\BigStrut
  $cD^{-1}(\log D)^{-2}$ & none & Masser (1989) \cite{MasserEC} \\ \hline
\BigStrut
  $cD^{-2/3}$ & $j$ nonintegral
       & Hindry-Silverman (1990) \cite{HindrySilverman} \\ \hline
\BigStrut
  $c$ & $j$ nonintegral or CM & Baker (2002) \cite{Baker} \\ \hline
\BigStrut
  $c$ & none & Silverman (2003) \\ \hline
\end{tabular}
\caption{History of lower bounds for $\hhat$ in $E(K^\ab)$}
\label{table:lowerbdinKabhistory}
\end{center}
\end{table}

\subsection{Lehmer's Conjecture for Abelian Varieties}
A number of authors have considered the analog of Lehmer's conjecture
for abelian varieties~$A/K$ of dimension~$g\ge2$. Let~$\Lcal$ be an
ample symmetric line bundle on~$A$, let~$P\in A(\Kbar)$ be a
nontorsion point, and let \text{$D=[K(P){:}K]$}. Masser~\cite{MasserAV} proves
that there is a constant $\k=\k(\dim A)$ so that
$$
  \hhat_{\Lcal}(P) \ge cD^{-\k}.
$$
If~$A$ has complex multiplication, then David and 
Hindry~\cite{DavidHindry} generalize Laurent's result to give 
a Dobrowoski-type estimate
$$
  \hhat(P) \ge cD^{-1}\bigl(\log\log(D)/\log(D)\bigr)^{\k}.
$$
Restricting to special types of fields, S. Zhang notes that the
equidistribution theorems proven in~\cite{Szpiro,Zhang} imply that
if~$L$ is a finite extension of a totally real field,
then~$A(L)_\tors$ is finite and $\hhat_\Lcal(P)\ge c$ for all
nontorsion~$P\in A(L)$. (See~\cite[Theorem~1.8]{Baker}.) In
particular, this is true for~$A(K^\ab)$ if~$K$ itself is totally
real. [The author would like to thank Michael Rosen for pointing out this
last fact.]

\section{Preliminaries}
\label{section:preliminaries}
In this section we set notation and recall some basic facts about
the decomposition of the canonical height into a sum of local heights.

\subsection{Notation}
\label{subsection:notation}
We set the following notation.

\begin{notation}
\item[$K/\QQ$]
a number field.
\item[$M_K$]
The set of absolute values on a number field~$K$  extending the usual
absolute values on~$\QQ$. We write $v(x)=-\log|x|_v$, and for
nonarchimedean~$v$, we denote the associated prime ideal by~$\gp_v$.
\item[$E/K$]
an elliptic curve defined over $K$.
\item[$\hhat$]
the canonical height $\hhat:E(\Kbar)\to\RR$ on~$E$, 
see~\cite[VIII~\S9]{SilvermanAEC}.
\item[$\lhat$]
the local canonical height (N\'eron function) 
$$
  \lhat:M_\Kbar\times E(\Kbar)\to\RR\cup\{\infty\},
$$ 
normalized as described in~\cite[VI.1.1]{SilvermanATAEC}.
\end{notation}

\begin{remark}
\label{remark:normalizations}
With absolute values normalized as above, the product formula reads
$$
  \sum_{v\in M_K} \frac{[K_v:\QQ_v]}{[K:\QQ]}v(\a)=0
  \qquad\text{for all $\a\in K^*$.}
$$
We also have the usual formula for finite extensions~$L/K$,
$$
  \sum_{w\in M_L,\,w|v} [L_w:K_v] = [L:K].
$$
\end{remark}

\begin{remark}
\label{remark:smallestvaluation}
If~$v$ is nonarchimedean, say corresponding to a prime
ideal~$\gp$ of residue characteristic~$p$ and ramification
index~$e=e(\gp/p)$, then~$v$ is discrete and its smallest postive
value is
$$
  \inf\{ v(\a) : \text{$\a\in K^*$ and $v(\a)>0$} \}
  = \frac{\log p}{e}.
$$
This is clear, since if~$\pi$ is a uniformizer in~$K$ for~$\gp$,
then locally~$(\pi)^e=(p)$, so $v(\pi)=v(p)/e=(\log p)/e$.
\end{remark}

\subsection{Local Height Functions}
\label{subsection:localhts}
We recall some well-known properties of the local
canonical height on an elliptic curve.

\begin{theorem}
\label{theorem:localheightproperties}
Let
$$
  \lhat:M_\Kbar\times E(\Kbar)\to\RR\cup\{\infty\},
$$
be the local canonical height (also called a N\'eron function) for the
divisor~$(O)$ and normalized as described
in~\cite[VI.1.1]{SilvermanATAEC}.
\begin{parts}
\Part{(a)}
For any finite extension $L/K$ and any point $P\in A(L)\setminus\{O\}$,
$$
  \hhat(P) =  \sum_{w\in M_L} \frac{[L_w:\QQ_w]}{[L:\QQ]} \lhat_{w}(P).
$$
\Part{(b)}
There is a constant $c_1(E/K)\ge0$ so that for any finite
extension~$L/K$ and any point $P\in A(L)\setminus\{O\}$,
$$
  \sum_{w\in M_L} \frac{[L_w:\QQ_w]}{[L:\QQ]} 
     \min\bigl\{\lhat_{w}(P),0\bigr\} \ge -c_1(E/K).
$$
\Part{(c)}
Let~$v\in M_K$ be a finite place of good reduction for~$E$, and choose
a minimal Weierstrass equation for~$E$ at~$v$,
$$
  y^2+a_1xy+a_3y=x^3+a_2x^2+a_4x+a_6.
$$
Let~$L/K$ be a finite extension and let~$w\in M_L$ lie over~$v$. Then
$$
  \lhat_{w}(P) = \max\{ w(x/y), 0 \}
  \qquad\text{for all $P=(x,y)\in E(L)\setminus\{O\}$.}
$$ 
\end{parts}
\end{theorem}

\begin{proof}
(a) 
This is~\cite[Theorem~VI.2.1]{SilvermanATAEC}.
\par
(c)
Since we are assuming that~$E$ has good reduction and that we have
chosen a minimal Weierstrass equation, this
is~\cite[Theorem~VI.4.1]{SilvermanATAEC}, except that we have 
written~$w(x/y)$ instead of~$\frac12w(x^{-1})$.
(Note that since~$E$ has good reduction and the
equation is minimal at~$v$, the same is true for any extension~$w$
of~$v$, and hence $w(\D)=0$ and~$E_0(L_w)=E(L_w)$.)
The integrality of the Weierstrass equation implies easily that
$$
  w(x^{-1})<0 \Longleftrightarrow w(x/y)<0
  \quad\text{and}\quad
  \min\bigl\{0,3w(x)\bigr\}=\min\bigl\{0,2w(y)\bigr\},
$$
so $\frac12w(x^{-1})=w(x/y)$ if either is negative.
\par
(b)
Let~$S_K\subset M_K$ be the union of the archimedean places of~$K$ and
the set of places of~$K$ at which~$E$ has bad reduction, and 
let~$S_L\subset M_L$ be the set of places of~$L$ lying above places of~$K$.
From~(c), we know that~$\lhat_w(P)\ge0$ for all $w\notin S_L$.
\par
Suppose now that $v\in S_K$ is archimedean. We fix an isomorphism
$E(\Kbar_v)\cong E(\CC)\cong\CC/q_v^\ZZ$ with $0<|q_v|_v<e^{-\pi}$.
Then the local
height is given by the explicit formula~\cite[VI.3.4]{SilvermanATAEC}
\begin{align}
  \label{equation:explicitformulaforlocalht}
  \lhat_v(u) =
  \frac12\BB_2&\left(\frac{\log|u|_v}{\log|q_v|_v}\right)
                 \log|q_v^{-1}|_v \\
   &-\log|1-u|_v-\sum_{n\ge1}\log|(1-q_v^nu)(1-q_v^n/u)|_v, \notag
\end{align}
where~$\BB_2(t)$ is the Bernoulli polynomial~$t^2-t+\frac16$ for~$0\le
t\le 1$, extended periodically to~$\RR/\ZZ$. It is easy to see
from~\eqref{equation:explicitformulaforlocalht} that
$$
  \lhat_v(u) \ge
  \frac12\BB_2\left(\frac{\log|u|_v}{\log|q_v|_v}\right)
                 \log|q_v^{-1}|_v - c_2
$$
for an absolute constant~$c_2$. The polynomial~$\BB_2(t)$ has a minimum at
$t=\frac12$, so we find that
$$
  \lhat_v(P) \ge -\frac{1}{24}\log|q_v^{-1}|_v - c_2
  \quad\text{for all $P\in E(\Kbar_v)\setminus \{O\}$.}
$$
Since~$K$ and~$E$ are fixed, this says that there is
a constant $c_3=c_3(E/K)$ so that $\lhat_v(P)\ge-c_3$ for all 
archimedean places~$v\in M_K$ and all points 
\text{$P\in E(\Kbar_v)\setminus \{O\}$}. It follows that for any finite
extension~$L/K$ and any point $P\in A(L)\setminus\{O\}$,
$$
  \sum_{w\in M_L^\infty} \frac{[L_w:\QQ_w]}{[L:\QQ]} 
     \min\bigl\{\lhat_{w}(P),0\bigr\} 
  \ge   \sum_{w\in M_L^\infty} \frac{[L_w:\QQ_w]}{[L:\QQ]}  (-c_3)
  = -c_3.
$$
Note that if $w$ lies over~$v$, then $\bar L_w\cong\Kbar_v$, so
we can compute the local height~$\lhat_w(P)$ using the inclusions
$E(L)\subset E(L_w)\subset E(\Kbar_v)$.
\par
Next suppose that~$v\in M_K$ is a nonarchimedean place at which~$E/K$
has split multiplicative reduction. Then the exact same argument gives
an analogous lower bound for~$\lhat_v$, because the local height is
given by the same formula~\eqref{equation:explicitformulaforlocalht},
see~\cite[VI.4.2]{SilvermanATAEC}.
\par
Finally, if~$E/K$ has additive or nonsplit multiplicative reduction 
at~$v$, then we may either make a finite extension of~$K$ to
reduce to one of the previous cases or use explicit formulas
in these cases (see, e.g.,~\cite[Table~1]{Silvermanthesis}) to obtain the
desired lower bound.
\end{proof}

\begin{remark}
Using explicit formulas for the local canonical height, it is not
difficult to give an explicit estimate for the constant~$c_1(E/K)$ in
terms of the $j$-invariant and minimal discriminant~$\Dcal_{E/K}$
of~$E/K$. Roughly, one can take~$c_1$ to have the form
$$
  c_1 = c_1'\max\bigl\{1,h(j),\log\Norm_{K/\QQ}\Dcal_{E/K}\bigr\},
$$
where~$c_1'$ is an absolute constant. See~\cite{SilvermanCompHt}
for a similar computation of explicit
constants associated to local height functions.
\end{remark}

\begin{remark}
We have chosen to prove Theorem~\ref{theorem:localheightproperties}(b)
by appealing to explicit formulas for the local height on elliptic
curves. It is also possible to prove this result by appealing to the
$M_K$-continuity of the local height and using an $M_K$-compactness
type of argument. Of course,~$E(\Kbar_v)$ is not even locally compact
for nonarchimedean places~$v$, so one must be careful. A detailed
development of the correct concept, which is called $M_K$-boundedness,
is given in~\cite{LangDG}.
\end{remark}

\section{Torsion Points in Abelian Extensions}
\label{section:torsion}
The following result is a simple consequence of a deep theorem of Serre.

\begin{theorem}[Serre \cite{SerreImageOfGalois}]
\label{theorem:serreimageofgalois}
Let $K/\QQ$ be a number field and let~$E/K$ be an elliptic curve
without~CM, that is, $\End(E/\Kbar)=\ZZ$. Then
$$
  \text{$E(K^\ab)_\tors$ is finite.}
$$
In particular,~$E(K^\ab)[\ell]=0$ for all but finitely many primes~$\ell$.
\end{theorem}
\begin{proof}
Serre~\cite{SerreImageOfGalois} proves that there is a finite set of 
primes~$S$ so that for all~$\ell\notin S$, the Galois representation
$$
  \rho_\ell:G_{\Kbar/K} \longrightarrow \Aut(E[\ell])\cong\GL_2(\FF_\ell)
$$
is surjective. Suppose that $T\in E(K^\ab)[\ell]$ with $T\ne O$ for
some $\ell\notin S$. The group~$\GL_2(\FF_\ell)$ acts transitively on
the nonzero vectors in~$\FF_\ell^2$, so Serre's theorem tells us that
the Galois orbit of~$T$ consists of all nonzero elements
of~$E[\ell]$. The conjugates of~$T$ are all defined over~$K^\ab$, so
we conclude that $E[\ell]\subset E(K^\ab)$. Hence the
represenation~$\rho_\ell$ factors through the maximal abelian
quotient~$G_{K^\ab/K}$ of~$G_{\Kbar/K}$. This contradicts the fact
that the image of~$\rho_\ell$ is a nonabelian group, which completes
the proof that~$E(K^\ab)[\ell]=0$ for all $\ell\notin S$.  This
statement suffices for our later applications. In order to prove the
stronger statement that~$E(K^\ab)_\tors$ is finite, it remains to show
that for any particular~$\ell$, the $\ell$-power torsion in~$E(K^\ab)$
is finite. This follows by a similar argument using the fact, also due
to Serre~\cite{Serre}, that the image~$G_{\Kbar/K}$ in
$\Aut(T_\ell(E))$ is open, and thus has finite index.
\end{proof}

\section{The Unramified Case}
\label{section:unramifiedcase}
In this section we prove the basic estimate required for the proof of
our main result in the case that the extension~$L/K$ is unramified at
a (small) prime~$\gp$ of~$K$. We begin with the observation that
there is an element of the group ring~$\ZZ[G_{\Kbar/K}]$ that
simultaneously annihilates~$E(K^\ab)$ modulo every prime of~$K^\ab$
lying above~$\gp$.

\begin{theorem}
\label{theorem:characpolyoffrob}
Let~$K/\QQ$ be a number field, let $\gp$ be a prime of~$K$,
and let~$E/K$ be an elliptic curve with good reduction at~$\gp$.
Let
$$
  \F_\gp(X) = \det\bigl(X-\Frob_\gp\bigm|T_\ell(E)\bigr)
   = X^2 - aX + q
$$
be the characteristic polynomial of Frobenius at~$\gp$. Thus
$\F_\gp(X)\in\ZZ[X]$ with $q=\Norm_{K/\QQ}\gp$ and $|a|\le2\sqrt{q}$
(see~\cite[Chapter~V]{SilvermanAEC}).
\par
Let $\bar\gp$ be a prime of~$\Kbar$ lying over~$\gp$ and
let~$\s\in(\bar\gp,L/K)\subset G_{\Kbar/K}$ be in the associated
Frobenius conjugacy class.
\begin{parts}
\Part{(a)}
For all $P\in E(\Kbar)$,
$$
  \F_\gp(\s)P \equiv O \pmod{\bar\gp}.
$$
(This congruence is taking place on the N\'eron model of~$E$, or
more prosaically, on a Weierstrass equation for~$E$ that is
minimal at~$\gp$.)
\Part{(b)}
If $P\in E(\Kbar)$ satisfies $\F_\gp(\s)P=O$, then~$P$ is a torsion point.
\end{parts}
\end{theorem}
\begin{proof}
(a)
When reduced modulo~$\bar\gp$, the element~$\s\in G_{\Kbar/K}$ acts as
the $q$-power Frobenius map $f_q\in\End(\tilde E/\FF_\gp)$. Further,
the map~$\F_\gp(f_q)$ annihilates~$T_\ell(\tilde E/\FF_\gp)$,
since~$\F_\gp$ is the characteristic polynomial of~$f_q$ acting
on~$T_\ell(\tilde E/\FF_\gp)$ and the Cayley-Hamilton theorem tells us
that a linear transformation satisfies its own characteristic equation.
However, we have the general fact that the map
$$
  \End(E) \longrightarrow \End(T_\ell(E))
$$
is injective (cf{.}~\cite[V.7.3]{SilvermanAEC}), so we conclude that
$\F_\gp(f_q)=0$ as an element of~$\End(\tilde E/\FF_\gp)$. In other
words,
$$
  \F_\gp(f_q)Q = O \qquad\text{for all $Q\in \tilde E(\bar\FF_\gp)$.}
$$
Finally, using the fact that the reduction map commutes with the action
of Galois, we see that for any~$P\in E(\Kbar)$, the point~$\F_\gp(\s)P$
is in the kernel of reduction. In other words,
$\F_\gP(\s)P\equiv O\MOD{\bar\gp}$, which completes the proof of~(a).
\par
(b)
Let $P\in E(\Kbar)$ satisfy $\F_\gp(\s)P=O$. Fix a finite Galois
extension~$L/K$ with~$P\in E(L)$, say of degree $m=[L{:}K]$.
Then~$\s^m=1$ in~$G_{L/K}$, so in particular, $\s^mP=P$. Let
$$
  r = \Resultant(\F_\gp(X),X^m-1)\in\ZZ.
$$
The complex roots of~$X^m-1$ have absolute value~$1$
and the complex roots of~$\F_\gp(X)$ have absolute value~$\sqrt{q}$,
so they have no complex roots in common. It follows that~$r\ne0$.
\par
The resultant of two polynomals~$\ZZ[X]$ is a generator for the ideal
that they generate, so we can
find~$a(X),b(X)\in\ZZ[X]$ satisfying
$$
  a(X)\F_\gp(X)+b(X)(X^m-1)=r.
$$
Substituting $X=\s$ gives the identity 
$$
  a(\s)\F_\gp(\s) + b(\s)(\s^m-1)=r
$$
in the group ring $\ZZ[G_{\Kbar/K}]$. Hence
$$
  rP = a(\s)\bigl(\F_\gp(\s)P\bigr)
      + b(\s)\bigl((\s^m-1)P\bigr) = O,
$$
so~$P$ is a point of finite order.
\end{proof}

\begin{corollary}
\label{corollary:unramifiedbound}
Continuing with the notation and assumptions from
Theorem~\ref{theorem:characpolyoffrob}, let $P\in E(\Kbar)$ be a nontorsion
point, and let~$e$ be the absolute ramification index of~$\bar\gp$ in
the field of definition~$K(P)$ of~$P$. (That is, let~$e$ be the ramification
index of the prime ideal~\text{$\bar\gp\cap K(P)$} over the prime~$p$.) Then
$$
  \lhat_{\bar\gp}(\F_\gp(\s)P) \ge \frac{\log p}{e}.
$$
\end{corollary}
\begin{proof}
To ease notation, let $Q=\F_\gp(\s)P$.
Theorem~\ref{theorem:characpolyoffrob}(a) tells us that $Q\equiv
O\MOD{\bar\gp}$, so~$Q$ is in the kernel of reduction
modulo~$\bar\gp$. Further, the assumption that~$P$ is nontorsion and
Theorem~\ref{theorem:characpolyoffrob}(b) imply that $Q\ne O$.
Let~$\bar v$ be the absolute value associated to~$\bar\gp$.
Then on a minimal Weierstrass equation,
we have $3{\bar v}(x(Q))=2{\bar v}(y(Q))<0$, so we can apply
Theorem~\ref{theorem:localheightproperties}(c) to conclude that
$$
  \l_{\bar\gp}(Q) = {\bar v}(x(Q)/y(Q)) > 0.
$$
Finally, we use the fact that for any extension~$L/K$, the minimum
positive value of~$\bar v$ restricted to~$L$ is $(\log p)/e$, where~$e$
is the ramification index of~$\bar\gp$ in~$L$. (See
Remark~\ref{remark:smallestvaluation}.)
\end{proof}

\section{The Ramified Case}
\label{section:ramifiedcase}
In this section we prove the basic estimate needed to handle the case
of ramified extensions. We begin by recalling the proof of the key
number field lemma, which is due to Amoroso and Dvornicich. We then
apply this lemma to obtain an analogous result for elliptic curves.
In~\cite{Baker}, this was done by altering the given point to make sure
it is nonzero modulo~$\gP$. We describe an alternative
approach in which we force the point to be zero modulo~$\gP$. This
allows us to work entirely within the formal group, where computations
are much easier.

\subsection{A congruence in ramified abelian extensions}

The following lemma is modeled
after~\cite[Lemma~2]{AmorosoDvornicich}.  See
also~\cite[Lemma~3.2]{AmorosoZannier} and \cite[Lemma~3.5]{Baker}.

\begin{lemma}[Amoroso-Dvornicich \cite{AmorosoDvornicich}]
\label{lemma:inertiacongruence}
Let $K/\QQ$ be a number field, let~$\gp$ be a degree~1 prime in~$K$
with residue characteristic~$p$, and let~$L/K$ be an abelian extension
that is ramified at~$\gp$. Let~$\gP$ be a prime of~$L$ lying
over~$\gp$, and let~$\Ocal_{L,\gP}$ denote teh localization of~$L$
at~$\gP$.  Then there exists an element~$\t\in I_{L/K}$ with~$\t\ne1$
such that
$$
  \t(\a)^p \equiv \a^p \pmod{p\Ocal_{L,\gP}}
  \qquad\text{for all $\a\in\Ocal_{L,\gP}$.}
$$
(Note that the strength of this result is that the congruence is
modulo~$p$, and not merely modulo~$\gP$.)
\end{lemma}
\begin{proof}
Without loss of generality, we may replace~$K$ and~$L$ by their
completions~$K_\gp=\QQ_p$ and~$L_\gP$, respectively. Then~$L$ is an
abelian extension of~$\QQ_p$, so by the local Kronecker-Weber theorem,
there is an integer~$m\ge1$ so
that~$L\subset\QQ_p(\z_m)$. (Here~$\z_m$ is a primitive
$m^{\text{th}}$~root of unity.) We take~$m$ to be minimal, i.e.,~$m$
is the conductor of~$L/\QQ_p$.
\par
The extension~$L/\QQ_p$ is ramified by assumption, which implies
that $p|m$. Let~$\t$ be a generator for the cyclic group
$G_{\QQ_p(\z_m)/\QQ_p(\z_{m/p})}$. The mimimality of~$m$ implies that
$L\not\subset\QQ_p(\z_{m/p})$, so the restriction of~$\t$ to~$G_{L/K}$
is not trivial. Further, since~$\t$ fixes~$\z_m^p$, we see that
$$
  \t(\z_m) = \o\z_m
  \qquad\text{for some primitive $p^{\text{th}}$ root of unity~$\o$.}
$$
\par
Now let~$\a\in\Ocal_L\subset\Ocal_{\QQ_p(\z_m)}=\ZZ_p[\z_m]$. Thus
$\a=f(\z_m)$ for some polynomial $f(X)\in\ZZ_p[X]$, and hence
$$
  \t\a=\t(f(\z_m))=f(\t\z_m)=f(\o\z_m).
$$
Now taking $p^{\text{th}}$~powers yields
\begin{align*}
  (\t\a)^p
  = f(\o\z_m)^p
 &\equiv f\bigl((\o\z_m)^p\bigr)   \pmod{p\ZZ_p[\z_m]} \\
 &= f(\z_m^p) \\
 &\equiv f(\z_m)^p   \pmod{p\ZZ_p[\z_m]}\\
 &= \a^p.
\end{align*}
Since $L\cap p\ZZ_p[\z_m]=p\Ocal_L$, this completes the proof.
\end{proof}

\subsection{Ramified points in $E(K^\ab)$}
In this section we prove the key estimate needed to handle the case of
ramified extensions. In particualr, by applying \text{$(\t-1)^2$} to a
point, we can do most of our computations in the formal group, 
where the multiplication-by-$p$ map is easier to describe.

\begin{lemma}
\label{lemma:ramifiedcase}
Let $E/K_\gp$ be an elliptic curve defined over a local field
$K_\gp/\QQ_p$ and assume that~$E$ has good reduction at~$\gp$.
Let~$L_\gP/K_\gp$ be a finite Galois extension that is ramified
at~$\gp$, let~$\Ocal_\gP$ be the ring of integers of~$L_\gP$, and
let~$\t\in I_{L_\gP/K_\gp}$ be in the inertia group
of~$L_\gP/K_\gp$. We denote by~$E_1(L_\gP)$ the kernel of reduction
and by~$\hat E$ the formal group of~$E$, so there is a natural
isomorphism $E_1(L_\gP)\cong\hat
E(\gP)$~\cite[Chapter~IV]{SilvermanAEC}.
\begin{parts}
\Part{(a)}
Let $P\in E(L_\gP)$. Then $(\t-1)P\in E_1(L_\gP)$.
\Part{(b)}
Fix a minimal Weierstrass equation for~$E$ and let~$z=-x/y$ be a
parameter for the formal group~$\hat E$. Then
\begin{align*}
  z\bigl([p](\t-1)Q\bigr) \in
  \bigl( (\t z(Q))^p - z(Q)^p \bigr)&\Ocal_\gP + p\Ocal_\gP \\
  &\quad\text{for all $Q\in E_1(L_\gP)\cong\hat E(\gP)$.}
\end{align*}
\end{parts}
\end{lemma}
\begin{proof}
(a) 
This is clear, since an element of the inertia group fixes everything 
modulo~$\gP$, so $\t P\equiv P\MOD{\gP}$, and hence the difference
$\t P-P$ lies in the kernel of reduction.
\par
(b) 
This is an immediate consequence of the following elementary
lemma about commutative formal groups, applied with $x=\t z(Q)$
and $y=z(Q)$.
\end{proof}

\begin{lemma}
\label{lemma:formalgroup}
Let~$R$ be a ring, let~$p\in\ZZ$ be a prime, and let
$$
  F(x,y) \in R[[x,y]]
$$
be a formal group over~$R$. Let~$\iota(t)\in R[[t]]$ be the inversion series 
for~$F$ and let~$M_p(t)\in R[[t]]$ be
the multiplication-by-$p$ series for~$F$. Then
$$
  M_p\bigl(F(x,\iota(y))\bigr) \in (x^p-y^p)R[[x,y]] + p R[[x,y]].
$$
\end{lemma}
\begin{proof}
There are power series $a(t),b(t)\in R[[t]]$ with $a(0)=b(0)=0$
so that multiplication-by-$p$ is given by power series of the form
$$
  M_p(t) = a(t^p) + pb(t).
$$
See~\cite[IV.4.4]{SilvermanAEC}. Further, the definition of~$\iota(t)$
implies that $F(x,\iota(y))$ vanishes when $y=x$, and hence
that~$F(x,\iota(y))$ is divisible by~$x-y$, say
$$
  F(x,\iota(y))=(x-y)G(x,y)\qquad\text{with}\qquad
  G(x,y)\in R[[x,y]].
$$
Then
\begin{align*}
  M_p\bigl(F(x,\iota(y))\bigr) 
  &\in a\bigl(F(x,\iota(y))^p\bigr) + p R[[x,y]] \\
  &\in a\bigl((x-y)^pG(x,y)^p\bigr) + p R[[x,y]] \\
  &\in (x-y)^p R[[x,y]]+ p R[[x,y]] \\
  &\in (x^p-y^p) R[[x,y]]+ p R[[x,y]] 
\end{align*}
\end{proof}

Combining Lemma~\ref{lemma:inertiacongruence} and
Lemma~\ref{lemma:ramifiedcase} gives the following crucial local
contribution to the canonical height in the ramified case. The key is
the fact that the lower bound does not depend on the ramification
degree of the field of definition of the point.

\begin{proposition}
\label{proposition:ramifiedcontribution}
Let $K/\QQ$ be a number field, let $P\in E(K^\ab)$, and let $L=K(P)$
be the field of definition of~$P$. Fix a degree~1 unramified
prime~$\gp$ of~$K$, let $p=\Norm_{K/\QQ}\gp$, and suppose that~$\gp$ ramifies
in~$L$. Then there exists an element~$\t\in I_{L/K}$ with~$\t\ne1$
such that the point
$$
  P' = [p]\bigl((\t-1)^2P\bigr)
$$
satisfies
$$
  \lhat_\gP(P')\ge\log p
$$
for all primes~$\gP$ of~$L$ lying over~$\gp$. (Note that if~$P'=O$, then
$\lhat_\gP(P')=\infty$ by definition, so the result is vacuously true.)
\end{proposition}
\begin{proof}
The extension~$L/K$ is abelian, so Lemma~\ref{lemma:inertiacongruence}
says that there is a nontrivial element $\t\in I_{L/K}$ such that
\begin{equation}
  \label{equation:modpcongruence}
  (\t\a)^p \equiv \a^p \pmod{p\Ocal_L}
  \qquad\text{for all $\a\in\Ocal_L$.}
\end{equation}
Let $Q=(\t-1)P$, so $P'=[p](\t-1)Q$. Lemma~\ref{lemma:ramifiedcase}(a)
says that~$Q$ is in the kernel of reduction modulo~$\gP$. Fix a
minimal Weierstrass equation for~$E$ at~$\gp$ and let~$z=-x/y$.  
Note that~$z(Q)$ is in the localization~$\Ocal_{L,\gP}$ of~$\Ocal_L$
at~$\gP$. Lemma~\ref{lemma:ramifiedcase} tells us that
$$
  z(P') \in \bigl((\t z(Q))^p-z(Q)^p\bigr)\Ocal_{L,\gP}+p\Ocal_{L,\gP}.
$$
Applying~\eqref{equation:modpcongruence} with~$\a=z(Q)$, we find
that~$z(P')\in p\Ocal_{L,\gP}$.
Hence Theorem~\ref{theorem:localheightproperties}(c) yields
$$
  \lhat_\gP(P') = v_\gP(z(P')) \ge v_\gP(p) = \log p.
$$
\end{proof}

\section{Proof of the Main Theorem}
\label{section:proofofmaintheorem}
We are now ready to prove our main result, which we restate for the
convenience of the reader.

\begin{theorem}
\label{theorem:maintheoremrestated}
Let $K/\QQ$ be a number field, let~$E/K$ be an elliptic curve, 
and let~$\hhat:E(\Kbar)\to\RR$ be the canonical height on~$E$.
There is a constant $C(E/K)>0$ such that every nontorsion point
$P\in E(K^\ab)$ satisfies
$$
  \hhat(P) > C.
$$
\end{theorem}
\begin{proof}
Baker~\cite{Baker} has proven
Theorem~\ref{theorem:maintheoremrestated} in the case that the
elliptic curve has complex multiplication, so we may assume
that~$\End(E/\Kbar)=\ZZ$. (Baker also proves the theorem in the case
that~$j(E)$ is nonintegral, but our proof will cover all non-CM curves.)
\par 
We begin by fixing a rational prime~$p$ and a prime~$\gp$ of~$K$ lying
over~$p$ with the following properties:
\begin{enumerate}
\item
$E(K^\ab)[p] = \{O\}$.
\item
$p\ge \exp\bigl([K:\QQ](1+c_1(E/K))\bigr)$, where $c_1(E/K)$ is the
constant appearing in Theorem~\ref{theorem:localheightproperties}(b).
\item
$E$ has good reduction at~$\gp$.
\item
$\gp$ is an unramified prime of degree~1.
\end{enumerate}
Our assumption that~$\End(E/\Kbar)=\ZZ$ and Serre's
Theorem~\ref{theorem:serreimageofgalois} tell us that~(1) eliminates
only finitely many primes~$p$, and similarly~(2) and~(3) eliminate only
finitely many primes~$\gp$. The condition~(4)
describes a set of primes of density
one, so we can find a prime satisfying~(1)--(4) that depends
only on the field~$K$ and the curve~$E/K$. Note,
however, that we are strongly using here the assumption that~$E$ does
not have~CM, since if~$E$ has~CM, then~(1) may be false for all primes.
\par
Let $P\in E(K^\ab)$ be a nontorsion point and let $L=K(P)$ be its field
of definition. The proof of the theorem proceeds by induction on the
ramification degree of~$L/K$ at~$\gp$.
\par
We begin with the unramified case, say $\gp\Ocal_L=\gP_1\gP_2\cdots\gP_d$.
The extension~$L/K$
is abelian, so the Frobenius elements associated to the~$\gP_i$
are all the same. Let $\s=(\gP_i,L/K)\in G_{L/K}$ be this Frobenius element,
and let~$\F_\gp(X)=X^2-aX+p$ be the characteristic polynomial 
of Frobenius. Then Theorem~\ref{theorem:characpolyoffrob}(b) tells
us that~$\F_\gp(\s)P$ is a nontorsion point, so
Corollary~\ref{corollary:unramifiedbound} and our
assumption that~$\gP_i$ is unramified  imply that
$$
  \lhat_{\gP_i}\bigl(\F_\gp(\s)P\bigr) \ge \log p.
$$
Note that this is true for every~$\gP_i$, because 
every~$\gP_i$ has the same~$\s\in
G_{L/K}$ as its Frobenius element. Adding over
the~$\gP_i$, we obtain the estimate
\begin{equation}
  \label{equation:unramifiedpcontribution}
  \sum_{\gP|\gp} \frac{[L_\gP:K_\gp]}{[L:K]}
       \lhat_{\gP}\bigl(\F_\gp(\s)P\bigr)
  \ge \log p.
\end{equation}
\par
Next we combine this positive lower bound with the (potentially negative)
contribution from the other absolute values. We compute:
\begin{align*}
  \hhat\bigl(&\F_\gp(\s)P\bigr) \\
   &= \sum_{w\in M_L} \frac{[L_w:\QQ_w]}{[L:\QQ]}\lhat_w\bigl(\F_\gp(\s)P\bigr)
    \quad\text{(Theorem~\ref{theorem:localheightproperties}(a))} 
    \allowdisplaybreaks\\
  & \ge
    \sum_{w\in M_L,\, w|v_\gp} 
            \frac{[L_w:\QQ_w]}{[L:\QQ]}\lhat_w\bigl(\F_\gp(\s)P\bigr) 
    \allowdisplaybreaks\\
  &\qquad\qquad\qquad\qquad {}+ \sum_{w\in M_L,\, w\notdivide v_\gp} 
            \frac{[L_w:\QQ_w]}{[L:\QQ]}
            \min\bigl\{\lhat_w\bigl(\F_\gp(\s)P\bigr),0\bigr\} 
    \allowdisplaybreaks\\
  & \ge
    \biggl(\sum_{w\in M_L,\, w|v_\gp} 
            \frac{[L_w:\QQ_w]}{[L:\QQ]}\lhat_w\bigl(\F_\gp(\s)P\bigr)\biggr)
    -c_1(E/K) 
    \quad\text{(Theorem~\ref{theorem:localheightproperties}(b))} 
    \allowdisplaybreaks\\
  & \ge  \frac{\log p}{[K:\QQ]} - c_1(E/K)
    \quad\text{(from \eqref{equation:unramifiedpcontribution})} \\
  & \ge 1
    \quad\text{from the choice of $p$.}
\end{align*}
\par 
In order to obtain a lower bound for~$\hhat(P)$, we use the fact
that~$\hhat$ is a Galois invariant~\cite[VIII.5.10]{SilvermanAEC}
positive semidefinite quadratic
form~\cite[VIII.9.3]{SilvermanAEC}. Hence
\begin{align*}
  \hhat\bigl(\F_\gp(\s)P\bigr)
  & =\hhat\bigl(\s^2P-[a]\s P+[p]P\bigr) \\
  & \le 3\bigl( \hhat(\s^2P) + \hhat([a]\s P) + \hhat([p]P) \bigr) \\
  &  = 3\bigl( \hhat(P) + a^2\hhat(P) + p^2\hhat(P) \bigr) \\
  & \le 3(1+4p+p^2)\hhat(P)\quad\text{since $|a|\le2\sqrt{p}$.}
\end{align*}
This gives the lower bound
$$
  \hhat(P) \ge \frac{1}{3(1+4p+p^2)},
$$
and since~$p$ was chosen depending only on~$E/K$, independent of the 
point~$P\in E(K^\ab)$, this completes the proof of 
Theorem~\ref{theorem:maintheoremrestated}
in the case that the extension~$K(P)/K$ is unramified at~$\gp$.
\par
Assume now that~$L/K$ is ramified at~$\gp$
and assume by induction that the theorem is
proven for all points defined over abelian extensions whose
$\gp$-ramification index is strictly smaller than~$e_\gp(L/K)$.
Let~$\t\in I_\gp(L/K)$ be the nontrivial element in the inertia
group described in Proposition~\ref{proposition:ramifiedcontribution}, 
so the point
$$
  P' = [p]\bigl((\t-1)^2P\bigr)
$$
satisfies~$\lhat_\gP(P')\ge\log p$ for all primes~$\gP|\gp$. Summing
over these primes, we find that
$$
  \sum_{\gP|\gp} \frac{[L_\gP:K_\gp]}{[L:K]} \lhat_\gP(P')
  \ge \log p.
$$
Note that this is exactly the same inequality that we obtained earlier
(cf.~\eqref{equation:unramifiedpcontribution}), except that
it applies to the point~$P'$, rather than to the point~$\F_\gp(\s)P$.
Hence as long as we know that~$P'\ne O$, then
the  same computation as was done above for~$\F_\gp(\s)P$
yields the same lower bound for~$P'$, namely
$$
  \hhat(P') \ge 1.
$$
In the other direction, we can estimate
\begin{align*}
  \hhat(P') = \hhat\bigl([p](\t-1)^2P\bigr)
  &= p^2\hhat(\t^2P-[2]\t P + P) \\
  &\le 3p^2\bigl(\hhat(\t^2P)+4\hhat(\t P)+\hhat(P)\bigr) \\
  &= 18p^2\hhat(P).
\end{align*}
Hence $\hhat(P)\ge 1/(18p^2)$, which completes the proof of
Theorem~\ref{theorem:maintheoremrestated} provided~$P'\ne O$.
\par
Finally, suppose that
$$
  P' = [p](\t-1)^2P = O.
$$
Let~$m$ be the order of~$\t$ in~$G_{L/K}$. There are
polynomials~$a(X),b(X)\in\ZZ[X]$ satisfying
$$
  a(X)(X^m-1) + b(X)p(X-1)^2 = mp(X-1).
$$
(The existence of such an identity follows immediately from the
fact that the resultant of~$X^{m-1}+\cdots+X+1$ and~$X-1$ is~$m$.
Explicitly, one can take~$a(X)=p$ and $b(X)=-\sum_{i=0}^{m-2}(m-1-i)X^i$.) 
We evaluate this identity at~$X=\t$ and apply it to the point~$P$.
We know that~$\t^m=1$ and $[p](\t-1)^2P=O$, so we find that
$$
  [mp](\t-1)P = O.
$$
In particular, $(\t-1)P\in E(L)_\tors$. However, we also 
know from Lemma~\ref{lemma:ramifiedcase} that \text{$(\t-1)P$} is
in the kernel of reduction modulo~$\gP$, so it follows
from general facts about formal groups~\cite[IV.3.2(b)]{SilvermanAEC}
that the order of~\text{$(\t-1)P$} is a power of~$p$. We now
recall that~$p$ was chosen so that~$E(K^\ab)[p]=\{O\}$, 
from which we conclude that \text{$(\t-1)P=O$}.
\par
Thus~$P$ is fixed by~$\t$, so~$P$ lies in the proper subfield~$L^\t$
of~$L$. (Note that~$\t\ne1$, so~$L^\t\ne L$.) Further,~$P$ is a
nontorsion point by assumption. Hence by induction we conclude
that~$\hhat(P)\ge C(E/K)$, which completes the proof of the theorem.
\end{proof}


\section{Generalization to Abelian Varieties}
\label{section:abelianvarieties}
Much of material in this paper can be generalized to the case of
abelian varieties, at the usual cost of making the arguments more
complicated. In this section we will sketch the changes required to
prove the following generalization of
Theorem~\ref{theorem:maintheoremrestated}.

\begin{theorem}
\label{theorem:maintheoremAV}
Let $K/\QQ$ be a number field, let~$A/K$ be an abelian variety,
let~$\Lcal$ be a symmetric ample line bundle on~$A/K$, and let
$$
  \hhat_\Lcal:A(\Kbar)\to\RR
$$
be the canonical height on~$A$ associated to~$\Lcal$.  Assume
that~$A/K$ contains no abelian subvarieties having complex
multiplication.  Then there is a constant $c=c(A/K,\Lcal)>0$ such that
for all nontorsion points $P\in A(K^\ab)$,
$$
  \hhat_\Lcal(P) \ge c.
$$
\end{theorem}
\begin{proof} (Sketch)
Replacing~$K$ by a finite extension, we may assume that~$A$ splits
into a product of geometrically simple abelian varieties, and then by
looking at the indivdual factors, we may assume that~$A$ itself is
geometrically simple. Next, replacing~$\Lcal$ by~$\Lcal^{\otimes n}$,
we may assume that~$\Lcal$ is very ample. Fix effective
divisors~$D_1,\ldots,D_r$ for~$\Lcal$ whose intersection consists of
only the point~$O$, and then fix canonical local heights (also known
at N\'eron functions)
$$
  \lhat_{D_i} : M_\Kbar\times A(\Kbar)\to\RR\cup\{\infty\}.
$$
See~\cite[Chapters~10 and~11]{LangDG} for standard properties of local
and global height functions. In particular, standard properties of
(canonical) local height functions give the analog of
Theorem~\ref{theorem:localheightproperties}.  Thus~\cite[Chapter~11,
Theorem~1.6]{LangDG} says that for any nonzero point~$Q\in A(\Kbar)$, we can
compute~$\hhat_\Lcal(Q)$ as 
$$
  \hhat_\Lcal(Q)=\sum_{w\in M_{K(Q)}}
  \frac{[K(Q)_v:\QQ_v]}{[K(Q):\QQ]}\lhat_{D_i}(w,Q)
$$
by choosing an index~$i$ with $Q\notin\Support(D_i)$, which gives the
analog of Theorem~\ref{theorem:localheightproperties}(a).  Similarly,
the analog of Theorem~\ref{theorem:localheightproperties}(b) follows
fairly directly from the $M_K$-positivity of local height functions
attached to positive divisors~\cite[Chapter~10,
Proposition~3.1]{LangDG}, and the analog of
Theorem~\ref{theorem:localheightproperties}(c) follows from the fact
that for places of good reduction, the canonical local height is given
as an intersection index on the special fiber of the N\'eron model
of~$A$~\cite[Chapter~11, Theorem~5.1]{LangDG}.
\par
Next consider Theorem~\ref{theorem:serreimageofgalois}, which uses a
deep result of Serre to limit the torsion in~$E(K^\ab)$.  For the
proof of Theorem~\ref{theorem:maintheoremAV}, we really only need to
know that
\begin{equation}
  \label{equation:noptorsionAV}
  A(K^\ab)[p] = 0\quad\text{for all but finitely many primes $p$.}
\end{equation}
(Indeed, it would suffice to know that $A(K^\ab)[p]=0$ for infinitely
many primes.) In certain cases, Serre has proven an ``image of
Galois'' theorem for abelian
varieties~\cite{SerreImageOfGaloisAV}. Fortunately, the weaker
statement~\eqref{equation:noptorsionAV} that we require is proven in
full generality by Zarhin~\cite{Zarhin} (see also~\cite{Ruppert}).
All of these results~\cite{Ruppert,SerreImageOfGaloisAV,Zarhin} rely
heavily on the  groundbreaking methods used by Faltings in his
proof of Tate's Isogeny Conjecture~\cite{Faltings}.
\par
Theorem~\ref{theorem:characpolyoffrob} generalizes directly to abelian
varieties. The characteristic polynomial~$\F_\gp(X)$ of~$\Frob_\gp$
acting on~$T_\ell(A)$ is a monic polynomial of degree~$2g$ whose roots
have norm~$\sqrt{q}$.  See~\cite[Section 18]{Milne}
or~\cite[Section~21, Application~II]{Mumford}.  The proof of both
parts of Theorem~\ref{theorem:characpolyoffrob} is an easy consequence
of these facts and the injectivity of
$$
  \End(\tilde A/\FF_{\bar\gp})\longrightarrow\End(T_\ell(\tilde A)).
$$
Further, if we write~~$\F_\gp(X)=\sum a_iX^{2g-i}$,
then the triangle inequality yields $|a_i|\le\binom{2g}{i}q^{i/2}$.
This allows us to handle the unramified case.
\par
The ramified case again relies on
the Amo\-ro\-so-Dvor\-ni\-cich Lemma~\ref{lemma:inertiacongruence},
together with
an abelian variety version of Lemma~\ref{lemma:ramifiedcase}.
And just as in the case of elliptic curves, 
Lemma~\ref{lemma:ramifiedcase} for abelian
varieties follows
from standard facts about commutative formal groups.
See~\cite{Frohlich} or~\cite{Hazewinkel} for basic results
on formal groups, and in particular a description of the 
multiplication-by-$p$ map that gives 
a higher dimensional version of Lemma~\ref{lemma:formalgroup}. 
\par
Having now assembled all of the required pieces, they fit together to
give a proof of Theorem~\ref{theorem:maintheoremAV} using the same
argument as given in the proof of
Theorem~\ref{theorem:maintheoremrestated}.
\end{proof}

\begin{remark}
Combining the methods of this paper, which essentially handles the
non-CM case, with the method for CM elliptic curves described in
Baker's article~\cite{Baker}, one may be able to prove in full
generality the inequality $\hhat_\Lcal(P)\ge C(A/K,\Lcal)$ for all
nontorsion points $P\in A(K^\ab)$ on all abelian varieties~$A/K$.
(This is essentially Conjecture~1.10 in~\cite{Baker}.)  The
details of this argument will be given in a subsequent
publication.
\end{remark}



\begin{thebibliography}{99}

\itemsep=\smallskipamount

\bibitem{AmorosoDvornicich}
F. Amoroso, R. Dvornicich, A lower bound for the height in abelian
extensions, \emph{J. Number Theory} \textbf{80} (2000), 260--272.

\bibitem{AmorosoZannier}
F. Amoroso, U. Zannier, A relative Dobrowolski lower bound over
abelian extensions, \emph{Ann. Scuola Norm. Sup. Pisa  Cl. Sci.} (4)
\textbf{29} (2000), 711--727.

\bibitem{AndersonMasser}
M. Anderson, D. Masser, Lower bound for heights on elliptic
curves, \emph{Math. Zeit.} \textbf{174} (1980), 23--34.


\bibitem{Baker}
M. Baker, Lower bounds for the canonical height on elliptic curves
over abelian extensions, preprint, December 2002.

\bibitem{BlanksbyMontgomery}
P.E. Blanksby, H.L. Montgomery, Algebraic integers near the unit
circle, \emph{Acta Arith.} \textbf{18} (1971), 355--369.

\bibitem{Dobrowolski} 
E. Dobrowolski, On a question of Lehmer and the number of irreducible
factors of a polynomial, \emph{Acta Arith.}  \textbf{34} (1979),
391--401.

\bibitem{DavidHindry}
S. David, M. Hindry,
Minoration de la hauteur de N\'eron-Tate sur les vari\'et\'es
ab\'eliennes de type C.M., \emph{J. Reine Angew. Math.} \textbf{529}
(2000), 1--74.


\bibitem{Faltings} 
G. Faltings, 
Endlichkeitss\"atze f\"ur abelsche Variet\"aten \"uber Zahlk\"orpern,
\emph{Invent. Math.} \textbf{73} (1983), 349--366.

\bibitem{Frohlich}
A. Fr\"ohlich, \emph{Formal Groups}, Lecture Notes in
Math. 74, Springer Verlag, New York 1968.


\bibitem{Hazewinkel}
M. Hazewinkel. \emph{Formal Groups and Applications}, 
Academic Press Inc., New York, 1978.

\bibitem{HindrySilverman}
M. Hindry and J.H. Silverman, On Lehmer's conjecture for elliptic
curves, in S\'eminaire de Th\'eorie des Nombres, Paris (1988-1989),
103--116, \emph{Progress in Math.} \textbf{91}, Birkh\"auser Boston, 
Boston, MA, 1990 103--116.


\bibitem{LangDG}
S. Lang, \emph{Fundamentals of Diophantine Geometry}, Springer-Verlag, 
New York, 1983.

\bibitem{Laurent}
M. Laurent, Minoration de la hauteur de N\'eron-Tate, \emph{S\'eminaire
de Th\'eorie de Nombres Paris 1981--1982}
\textbf{38} (1983), Birkh\"auser, Boston-Basil-Stuttgart, 137--152.

\bibitem{Lehmer}
D.H. Lehmer, Factorization of certain cyclotomic functions, 
\emph{Ann. of Math.} \textbf{34} (1933), 461--479. 

\bibitem{LehmerConjHomePage}
Lehmer Conjecture Website,
\texttt{http://www.math.ucla.edu/\char`\~mjm/lc/lc.html}

\bibitem{MasserEC}
D. Masser, Counting points of small height on elliptic curves,
\emph{Bull. Soc. Math. France} \textbf{117} (1989), 247--265.

\bibitem{MasserAV}
\bysame
Small values of the quadratic part of the N\'eron-Tate height on
an abelian variety, \emph{Compositio Math.} \textbf{53} (1984), 153--170.

\bibitem{Milne}
J.S. Milne, Abelian varieties, in \emph{Arithmetic Geometry}, ed. G. Cornell
and J.H. Silverman, Springer-Verlag, Berlin and New York, 1986.

\bibitem{Mumford}
D. Mumford, \emph{Abelian Varieties}, Tata Institute of Fundamental
Research Studies in Mathematics, No. 5,  Tata Institute, Bombay, 1970.

\bibitem{Pacheco}
A. Pacheco, Lehmer's conjecture in positive characteristic, 
preprint, April 2002.

\bibitem{Ruppert}
W.M. Ruppert,
Torsion points of abelian varieties in abelian extensions, 
unpublished manuscript,\newline
\texttt{http://www.math.uiuc.edu/Algebraic-Number-Theory/0101/tpavae.dvi}

\bibitem{SerreImageOfGalois}
J.-P. Serre, Propri\'et\'es galoisiennes des points d'ordre fini des
courbes elliptiques, \emph{Invent. Math.} \textbf{15} (1972),
259--331.

\bibitem{Serre}
\bysame
\emph{Abelian $l$-adic Representations and Elliptic Curves},
(revised reprint of the 1968 original), Research Notes in Mathematics, 
7, A K Peters, Ltd., Wellesley, MA, 1998.


\bibitem{SerreImageOfGaloisAV}
\bysame
\emph{Oeuvres Colllected Papers}, Volume IV, 1985--1998, 
136. R\'esum\'e des cours de 1985--1986 and 137. Lettre \'a Marie-France
Vign\'eras du 10/2/1986, Springer-Verlag, Berlin, 2000.


\bibitem{SilvermanAEC}
J.H. Silverman, \emph{The Arithmetic of Elliptic Curves},
Graduate Texts in Mathematics~106, Springer-Verlag, New York, 1986.

\bibitem{SilvermanATAEC}
\bysame
\emph{Advanced Topics in the Arithmetic of Elliptic Curves},
Graduate Texts in Mathematics~151, Springer-Verlag, New York, 1994.

\bibitem{Silvermanthesis}
\bysame
The N\'eron-Tate Height on Elliptic Curves, Ph.D. thesis, Harvard,
December 1981.

\bibitem{SilvermanDuke} 
\bysame
Lower bound for the canonical height on elliptic
curves, \emph{Duke Math. J.} \textbf{48} (1981), 633--648.

\bibitem{SilvermanCompHt}
\bysame
The difference between the Weil height and the canonical height on
elliptic curves, \emph{Math. Comp.} \textbf{55} (1990), 723--743.

\bibitem{Smyth}
C.J. Smyth, On the product of the conjugates outside the unit circle
of an algebraic number, \emph{Bull. London Math. Soc.} \textbf{3}
(1971), 169--175.

\bibitem{Stewart} 
C.L. Stewart, Algebraic integers whose conjugates lie near the unit
circle, \emph{Bull. Soc.  Math. France} \textbf{106} (1978), no. 2,
169--176.

\bibitem{Szpiro}
L. Szpiro, E. Ullmo, S. Zhang, \'Equir\'epartition des petis points,
\emph{Invent. Math.} \textbf{127} (1997), 337--347.

\bibitem{Zarhin}
Yu.G. Zarhin, 
Endomorphisms and torsion of abelian varieties.
\emph{Duke Math. J.} \textbf{54} (1987), 131--145.

\bibitem{Zhang}
S. Zhang, Equidistribution of small points on abelian varieties,
\emph{Annals of Math.} \textbf{147} (1998), 159--165.

\end{thebibliography}
\end{document}